\documentclass[12pt]{article}%
\usepackage{amsmath}
\usepackage{amsfonts}
\usepackage{amssymb}
\usepackage{graphicx}%
\providecommand{\U}[1]{\protect\rule{.1in}{.1in}}
\newtheorem{theorem}{Theorem}[section]

\newtheorem{conjecture}{Conjecture}[section]

\newtheorem{question}{Question}[section]

\providecommand{\boksie}{\ensuremath{\mathbin{\raisebox{0.3mm}{$\scriptstyle\square$}}}}
\begin{document}

\title{\textbf{Reconfiguration of Colourings and}\\\textbf{Dominating Sets in Graphs:\ a Survey}}
\author{C.M. Mynhardt\thanks{Supported by the Natural Sciences and Engineering
Research Council of Canada.}\\Department of Mathematics and Statistics\\University of Victoria, Victoria, BC, \textsc{Canada}\\{\small kieka@uvic.ca}
\and S. Nasserasr$^{\ast}$\\Department of Mathematics and Computer Science\\Brandon University, Brandon, MB, \textsc{Canada}\\{\small shahla.nasserasr@gmail.com}}
\date{}
\maketitle

\begin{abstract}
We survey results concerning reconfigurations of colourings and dominating
sets in graphs. The vertices of the $k$-colouring graph $\mathcal{C}_{k}(G)$
of a graph $G$ correspond to the proper $k$-colourings of a graph $G$, with
two $k$-colourings being adjacent whenever they differ in the colour of
exactly one vertex. Similarly, the vertices of the $k$-edge-colouring graph
$\mathcal{EC}_{k}(G)$ of $g$ are the proper $k$-edge-colourings of $G$, where
two $k$-edge-colourings are adjacent if one can be obtained from the other by
switching two colours along an edge-Kempe chain, i.e., a maximal two-coloured
alternating path or cycle of edges. 

The vertices of the $k$-dominating graph $\mathcal{D}_{k}(G)$ are the (not
necessarily minimal) dominating sets of $G$ of cardinality $k$ or less, two
dominating sets being adjacent in $\mathcal{D}_{k}(G)$ if one can be obtained
from the other by adding or deleting one vertex. On the other hand, when we
restrict the dominating sets to be minimum dominating sets, for example, we
obtain different types of domination reconfiguration graphs, depending on
whether vertices are exchanged along edges or not.

We consider these and related types of colouring and domination
reconfiguration graphs. Conjectures, questions and open problems are stated
within the relevant sections.

\end{abstract}

\section{Introduction}

\label{26CM-Sec_Intro} In graph theory, reconfiguration is concerned with
relationships among solutions to a given problem for a specific graph. The
reconfiguration of one solution into another occurs via a sequence of steps,
defined according to a predetermined rule, such that each step produces an
intermediate solution to the problem. The solutions form the vertex set of the
associated reconfiguration graph, two vertices being adjacent if one solution
can be obtained from the other in a single step. Exact counting of
combinatorial structures is seldom possible in polynomial time. Approximate
counting of the structures, however, may be possible. When the reconfiguration
graph associated with a specific structure is connected, Markov chain
simulation can be used to achieve approximate counting. Typical questions
about the reconfiguration graph therefore concern its structure
(connectedness\footnote{We use the term \emph{connectedness} instead of
\emph{connectivity} when referring to the question of whether a graph is
connected or not, as the latter term refers to a specific graph parameter.},
Hamiltonicity, diameter, planarity), realisability (which graphs can be
realised as a specific type of reconfiguration graph), and algorithmic
properties (finding shortest paths between solutions quickly).

Reconfiguration graphs can, for example, be used to study combinatorial Gray
codes. The term \textquotedblleft combinatorial Gray code\textquotedblright%
\ refers to a list of combinatorial objects so that successive objects differ
in some prescribed minimal way. It generalises Gray codes, which are lists of
fixed length binary strings such that successive strings differ by exactly one
bit. Since the vertices of a reconfiguration graph are combinatorial objects,
with two vertices being adjacent whenever they differ in some small way, a
Hamilton path in a reconfiguration graph corresponds to a combinatorial Gray
code in the source graph, and a Hamilton cycle to a cyclic combinatorial Gray code.

We restrict our attention to reconfigurations of graph colourings and
dominating sets (of several types). Unless stated otherwise, we use $n$ to
denote the order of our graphs. As is standard practice we denote the
chromatic number of a graph $G$ by $\chi(G)$, its clique number by $\omega
(G)$, and its minimum and maximum degrees by $\delta(G)$ and $\Delta(G)$,
respectively. We use $\gamma(G)$ and $\Gamma(G)$ to denote the domination and
upper domination numbers of $G$, that is, the cardinality of a minimum
dominating set and a maximum minimal dominating set, respectively.

One of the best studied reconfiguration graphs is the $k$-\emph{colouring
graph} $\mathcal{C}_{k}(G)$, whose vertices correspond to the proper
$k$-colourings of a graph $G$, with two $k$-colourings being adjacent whenever
they differ in the colour of exactly one vertex. When $\mathcal{C}_{k}(G)$ is
connected, a Markov process can be defined on it that leads to an
approximation of the number of $k$-colourings of $G$; this relationship
motivated the study of the connectedness of $\mathcal{C}_{k}(G)$. Some authors
consider list colourings with the same adjacency condition, while others
consider proper $k$-edge-colourings, where two $k$-edge-colourings of $G$ are
adjacent in the $k$-\emph{edge-colouring graph} $\mathcal{EC}_{k}(G)$ if one
can be obtained from the other by switching two colours along an
\emph{edge-Kempe chain}, i.e., a maximal two-coloured alternating path or
cycle of edges.

The domination reconfiguration graph whose definition most resembles that of
the $k$-colouring graph is the $k$-\emph{dominating graph} $\mathcal{D}%
_{k}(G)$, whose vertices are the (not necessarily minimal) dominating sets of
$G$ of cardinality $k$ or less, where two dominating sets are adjacent in
$\mathcal{D}_{k}(G)$ if one can be obtained from the other by adding or
deleting one vertex. The $k$-\emph{total-dominating graph} $\mathcal{D}%
_{k}^{t}(G)$ is defined similarly using total-dominating sets.

Other types of domination reconfiguration graphs are defined using only sets
of cardinalities equal to a given domination parameter $\pi$. For example, if
$\pi$ is the domination number $\gamma$, then the vertex set of the associated
reconfiguration graph, called the $\gamma$-\emph{graph} of $G$, consists of
the minimum dominating sets of $G$. There are two types of $\gamma$-graphs:
$\mathcal{J}(G,\gamma)$ and $\mathcal{S}(G,\gamma)$. In $\mathcal{J}%
(G,\gamma)$, two minimum dominating sets $D_{1}$ and $D_{2}$ are adjacent if
and only if there exist vertices $x\in D_{1}$ and $y\in D_{2}$ such that
$D_{1}-\{x\}=D_{2}-\{y\}$. The $\gamma$-graph $\mathcal{J}(G,\gamma)$ is
referred to as the $\gamma$-\emph{graph in the single vertex replacement
adjacency model} or simply the \emph{jump }$\gamma$\emph{-graph}. In
$\mathcal{S}(G,\gamma)$, two minimum dominating sets $D_{1}$ and $D_{2}$ are
adjacent if and only if there exist adjacent vertices $x\in D_{1}$ and $y\in
D_{2}$ such that $D_{1}-\{x\}=D_{2}-\{y\}$. The $\gamma$-graph $\mathcal{S}%
(G,\gamma)$ is referred to as the $\gamma$-\emph{graph in the slide adjacency
model} or the \emph{slide }$\gamma$\emph{-graph}. Note that $\mathcal{S}%
(G,\gamma)$ is a spanning subgraph of $\mathcal{J}(G,\gamma)$. In general we
define the \emph{slide }$\pi$\emph{-graph} similar to the slide $\gamma$-graph
and denote it by $\mathcal{S}(G,\pi)$.

We refer the reader to the well-known books \cite{26CM-CL} and
\cite{26CM-West} for graph theory concepts not defined here. Lesser known
concepts are defined where needed. We only briefly mention algorithmic and
complexity results, since a recent and extensive survey of this aspect of
reconfiguration is given by Nishimura \cite{26CM-Nishimura}. We state open
problems and conjectures throughout the text where appropriate.

\section{Complexity}

\label{SecComplex}Many of the published papers on reconfiguration problems
address complexity and algorithmic questions. The main focus of much of this
work has been to determine the existence of paths between different solutions,
that is, to determine which solutions are in the same component of the
reconfiguration graph, and if so, how to find a shortest path between two
solutions. The questions, therefore, are whether one solution is
\emph{reachable} from another according to the rules of adjacency, and if so,
to determine or bound the \emph{distance} between them. If all solutions are
reachable from one another, the reconfiguration graph is connected and its
diameter gives an upper bound on the distance between two solutions.

Complexity results concerning the connectedness and diameter of the
$k$-colouring graph $\mathcal{C}_{k}(G)$ are given in \cite[Section
6]{26CM-Nishimura}, and those pertaining to domination graphs can be found in
\cite[Section 7]{26CM-Nishimura}. We mention complexity results for
homomorphism reconfiguration in Section~\ref{SecWrochna}.

An aspect that has received considerable attention, but has not been fully
resolved, is to determine dividing lines between tractable and intractable
instances for reachability. Cereceda, Van den Heuvel, and Johnson
\cite{26CM-CHJ2} showed that the problem of recognizing bipartite graphs $G$
such that $\mathcal{C}_{3}(G)$ is connected is coNP-complete, but polynomial
when restricted to planar graphs. In \cite{26CM-CHJ3} they showed that for a
$3$-colourable graph $G$ of order $n$, both reachability and the distance
between given colourings can be solved in polynomial time. Bonsma and Cereceda
\cite{26CM-BC} showed that when $k\geq4$, the reachability problem is
PSPACE-complete. Indeed, it remains PSPACE-complete for bipartite graphs when
$k\geq4$, for planar graphs when $4\leq k\leq6$, and for bipartite planar
graphs when $k=4$. Moreover, for any integer $k\geq4$ there exists a family of
graphs $G_{N,k}$ of order~$N$ such that some component of $\mathcal{C}%
_{k}(G_{N,k})$ has diameter $\Omega(2^{N})$. Bonsma, Mouawad, Nishimura, and
Raman \cite{26CM-BMNR} showed that when $k\geq4$, reachability is strongly
NP-hard. Bonsma and Mouawad \cite{26CM-BM} explored how the complexity of
deciding whether $\mathcal{C}_{k}(G)$ contains a path of length at most $\ell$
between two given $k$-colourings of $G$ depends on $k$ and $\ell$, neatly
summarizing their results in a table. Other work on the complexity of
colouring reconfiguration include \cite{26CM-BB, 26CM-BJLPP, 26CM-BCHJ,
26CM-BMMN, 26CM-CHJ1, 26CM-DFFV, 26CM-FJP, 26CM-Ito2, 26CM-IKD, 26CM-Jerrum,
26CM-JKKPP, 26CM-JKKPP2, 26CM-Molloy}.

Haddadan, Ito, Mouawad, Nishimura, Ono, Suzuki, and Tebbal \cite{26CM-HIMNOST}
showed that determining whether $\mathcal{D}_{k}(G)$ is connected is
PSPACE-complete even for graphs of bounded bandwidth, split graphs, planar
graphs, and bipartite graphs, and they developed linear-time algorithms for
cographs, trees, and interval graphs. Lokshtanov, Mouawad, Panolan, Ramanujan,
and Saurabh \cite{26CM-LMPRS} showed that, although \texttt{W[1]}-hard when
parameterized by $k$, the problem is fixed-parameter tractable when
parameterized by $k+d$ for $K_{d,d}$-free graphs. For other works in this area
see \cite{26CM-MNRSS, 26CM-Tebbal}.

\section{Reconfiguration of Colourings}

\label{SecColour}The set of proper $k$-colourings of a graph $G$ has been
studied extensively via, for example, the Glauber dynamics Markov chain for
$k$-colourings; see e.g. \cite{26CM-DFFV, 26CM-DGM, 26CM-Jerrum, 26CM-LM,
26CM-Molloy}. Algorithms for random sampling of $k$-colourings and
approximating the number of $k$-colourings arise from these Markov chains. The
connectedness of the $k$-colouring graph is a necessary condition for such a
Markov chain to be rapidly mixing, that is, for the number of steps required
for the Markov chain to approach its steady state distribution to be at most a
polynomial in $\log(n)$, where $n=|V(G)|$.

\subsection{The $k$-Colouring Graph}

\label{Sec_k_colour}

Motivated by the Markov chain connection, a graph $G$ is said to be
$k$\emph{-mixing} if $\mathcal{C}_{k}(G)$ is connected. The minimum integer
$m_{0}(G)$ such that $G$ is $k$-mixing whenever $k\geq m_{0}(G)$ is called the
\emph{mixing number} of $G$. A $k$-colouring of $G$ is \emph{frozen} if each
vertex of $G$ is adjacent to at least one vertex of every other colour; a
frozen $k$-colouring is an isolated vertex of $\mathcal{C}_{k}(G)$. The
\emph{colouring number} $\operatorname{col}(G)$ of $G$ is the least integer
$d$ such that the vertices of $G$ can be ordered as $v_{1}\prec\cdots\prec
v_{n}$ so that $|\{v_{i}:i<j$ and $v_{i}v_{j}\in E(G)\}|<d$ for all
$j=1,...,n$. By colouring the vertices $v_{1},...,v_{n}$ greedily, in this
order, with the first available colour from $\{1,...,d\}$, we obtain a
$d$-colouring of $G$; hence $\chi(G)\leq\operatorname{col}(G)$. Here we should
mention that some authors define the colouring number to be $\max_{H\subseteq
G}\delta(H)$ where the maximum is taken over all subgraphs $H$ of $G$; this
number in fact equals $\operatorname{col}(G)-1$. Indeed, $\max_{H\subseteq
G}\delta(H)$ is often called the \emph{degeneracy} of $G$.


The choice of $k$ is important when we consider the connectedness and diameter
of $\mathcal{C}_{k}(G)$. Given two colourings $c_{1}$ and $c_{2}$, when $k$ is
sufficiently large each vertex can be recoloured with a colour not appearing
in either $c_{1}$ or $c_{2}$ and then recoloured to its target colour. Then
$\mathcal{C}_{k}(G)$ is connected and has diameter linear in the order of $G$.
This also shows that $m_{0}(G)$ is defined for each graph $G$. On the other
hand, if $k=2$ and $G$ is an even cycle, then no vertex can be recoloured and
$\mathcal{C}_{2}(G)=2K_{1}$.

Jerrum \cite{26CM-Jerrum} showed that $m_{0}(G)\leq\Delta(G)+2$ for each graph
$G$. Cereceda et al. \cite{26CM-CHJ1} used the colouring number to bound
$m_{0}$. Since $\operatorname{col}(G)\leq\Delta(G)+1$ and the difference can
be arbitrary, their result offers an improvement on Jerrum's
bound.\footnote{Bonsma and Cereceda \cite{26CM-BC} and Cereceda et al.
\cite{26CM-CHJ1} use the alternative definition of $\operatorname{col}(G)$; we
have adjusted their statements to conform to the definition given here.}

\begin{theorem}
\label{Thm_colG+1}\emph{\cite{26CM-CHJ1}}\hspace{0.1in}For any graph $G$,
$m_{0}(G)\leq\operatorname{col}(G)+1$.
\end{theorem}

Cereceda et al. \cite{26CM-CHJ1} used the graph $L_{m}=K_{m,m}-mK_{2}$ (the
graph obtained from the complete bipartite graph $K_{m,m}$ by deleting a
perfect matching) to obtain a graph $G$ and integers $k_{1}<k_{2}$ such that
$G$ is $k_{1}$-mixing but not $k_{2}$-mixing: colour the vertices in each
partite set of $L_{m}$ with the colours $1,...,m$, where vertices in different
parts that are ends of the same deleted edge receive the same colour. This
$m$-colouring is an isolated vertex in the $m$-colour graph $\mathcal{C}%
_{m}(L_{m})$. Hence $L_{m}$ is not $m$-mixing (there are many $m$-colourings
of $L_{m}$). They showed that for $m\geq3$, the bipartite graph $L_{m}$ is
$k$-mixing for $3\leq k\leq m-1$ and $k\geq m+1$ but not $k$-mixing for $k=m$.
They also showed that there is no expression $\varphi(\chi)$ in terms of the
chromatic number $\chi$ such that for all graphs $G$ and integers
$k\geq\varphi(\chi(G))$, $G$ is $k$-mixing.

Cereceda et al. \cite{26CM-CHJ1} also showed that if $\chi(G)\in\{2,3\}$, then
$G$ is not $\chi(G)$-mixing, and that $C_{4}$ is the only $3$-mixing cycle. In
contrast, for $m\geq4$ they obtained an $m$-chromatic graph $H_{m}$ that is
$k$-mixing whenever $k\geq m$: let $H_{m}$ be the graph obtained from two
copies of $K_{m-1}$ with vertex sets $\{v_{1},...,v_{m-1}\}$ and
$\{w_{1},...,w_{m-1}\}$ by adding a new vertex $u$ and the edges $v_{1}w_{1}$
and $\{uv_{i},uw_{i}\colon\,2\leq i\leq m-1\}$. In \cite{26CM-CHJ2}, the same
authors characterised $3$-mixing connected bipartite graphs as those that are
not foldable to $C_{6}$. [If $v$ and $w$ are vertices of a bipartite graph $G$
at distance two, then a \emph{fold} on $v$ and $w$ is the identification of
$v$ and $w$ (remove any resulting multiple edges); $G$ is \emph{foldable} to
$H$ if there exists a sequence of folds that transforms $G$ into $H$.]

Bonamy and Bousquet \cite{26CM-BB} used the Grundy number of $G$ to improve
Jerrum's bound on $m_{0}(G)$. A proper $k$-colouring of $G$ in colours
$1,...,k$ is called a \emph{Grundy colouring} if, for $1\leq i\leq k$, every
vertex with colour $i$ is adjacent to vertices of all colours less than $i$.
The \emph{Grundy number} $\chi_{g}(G)$ of a graph $G$ is the maximum number of
colours among all Grundy colourings of $G$. Note that $\chi_{g}(G)\leq
\Delta(G)+1$ and, as in the case of $\operatorname{col}(G)$, it can be
arbitrarily smaller.

\begin{theorem}
\label{ThmGrundy}\emph{\cite{26CM-BB}}\hspace{0.1in}For any graph $G$ of order
$n$ and any $k$ with $k\geq\chi_{g}(G)+1$, $\mathcal{C}_{k}(G)$ is connected
and $\operatorname{diam}(\mathcal{C}_{k}(G))\leq4n\chi(G)$.
\end{theorem}

Since the Grundy number of a \emph{cograph} (a $P_{4}$-free graph) equals its
chromatic number, Theorem \ref{ThmGrundy} implies that for $k\geq\chi(G)+1$, a
cograph $G$ is $k$-mixing and the diameter of $\mathcal{C}_{k}(G)$ is
$O(\chi(G)\cdot n)$, (i.e., linear in $n$). This result does not generalize to
$P_{r}$-free graphs for $r\geq5$. Bonamy and Bousquet constructed a family of
$P_{5}$-free graphs $\{G_{k}:k\geq3\}$ having both a proper $(k+1)$-colouring
and a frozen $2k$-colouring. They also showed that the graphs $L_{m}$
mentioned above are $P_{6}$-free with arbitrary large mixing number and asked
the following question.

\begin{question}
\emph{\cite{26CM-BB}}\hspace{0.1in}Given $r,k\in\mathbb{N}$, does there exist
$c_{r,k}$ such that for any $P_{r}$-free graph $G$ of order $n$ that is
$k$-mixing, the diameter of $\mathcal{C}_{k}(G)$ is at most $c_{r,k}\cdot n$?
\end{question}

Several other authors also considered the diameter of $\mathcal{C}_{k}(G)$ or
of its components when it is disconnected. Cereceda et al. \cite{26CM-CHJ3}
showed that if $G$ is a $3$-colourable graph with $n$ vertices, then the
diameter of any component of $\mathcal{C}_{3}(G)$ is $O(n^{2})$. In contrast,
for $k\geq4$, Bonsma and Cereceda \cite{26CM-BC} obtained graphs (which may be
taken to be bipartite, or planar when $4\leq k\leq6$, or planar and bipartite
when $k=4$) having $k$-colourings such that the distance between them is
superpolynomial in the order and size of the graph. They also showed that if
$G$ is a graph of order $n$ and $k\geq2\operatorname{col}(G)-1$, then
$\operatorname{diam}(\mathcal{C}_{k}(G))=O(n^{2})$. They stated the following conjecture.

\begin{conjecture}
\emph{\cite{26CM-BC}}\hspace{0.1in}For a graph $G$ of order $n$ and
$k\geq\operatorname{col}(G)+1$, $\operatorname{diam}(\mathcal{C}%
_{k}(G))=O(n^{3})$.
\end{conjecture}

Bonamy, Johnson, Lignos, Patel, and Paulusma \cite{26CM-BJLPP} determined
sufficient conditions for $\mathcal{C}_{k}(G)$ to have a diameter quadratic in
the order of $G$. They showed that $k$-colourable chordal graphs and chordal
bipartite graphs satisfy these conditions and hence have an $\ell$-colour
diameter that is quadratic in $k$ for $\ell\geq k+1$ and $\ell=3$,
respectively. Bonamy and Bousquet \cite{26CM-BB} proved a similar result for
graphs of bounded treewidth. Beier, Fierson, Haas, Russell, and Shavo
\cite{26CM-Beier} considered the girth $g(\mathcal{C}_{k}(G))$.

\begin{theorem}
\emph{\cite{26CM-Beier}}\hspace{0.1in}If $k>\chi(G)$, then $g(\mathcal{C}%
_{k}(G))\in\{3,4,6\}$. In particular, for $k>2$, $g(\mathcal{C}_{k}%
(K_{k-1}))=6$. If $k>\chi(G)+1$, or $k=\chi(G)+1$ and $\mathcal{C}_{k-1}(G)$
has an edge, then $g(\mathcal{C}_{k}(G))=3$. If $k=\chi(G)+1$ and $G\neq
K_{k-1}$, then $g(\mathcal{C}_{k}(G))\leq4$.
\end{theorem}


The Hamiltonicity of $\mathcal{C}_{k}(G)$ was first considered by Choo
\cite{26CM-ChooThesis} in 2002 (also see Choo and MacGillivray
\cite{26CM-ChooMacG}). Choo showed that, given a graph $G$, there is a number
$k_{0}(G)$ such that $\mathcal{C}_{k}(G)$ is Hamiltonian whenever $k\geq
k_{0}(G)$. The number $k_{0}(G)$ is referred to as the \emph{Gray code number}
of $G$, since a Hamilton cycle in $\mathcal{C}_{k}(G)$ is a (cyclic)
combinatorial Gray code for the $k$-colourings. Clearly, $k_{0}(G)\geq
m_{0}(G)$. By Theorem \ref{Thm_colG+1}, $m_{0}(G)\leq\operatorname{col}(G)+1$.
Choo and MacGillivray showed that one additional colour suffices to ensure
that $\mathcal{C}_{k}(G)$ is Hamiltonian.

\begin{theorem}
\emph{\cite{26CM-ChooMacG}}\hspace{0.1in}For any graph $G$ and $k\geq
\operatorname{col}(G)+2$, $\mathcal{C}_{k}(G)$ is Hamiltonian.
\end{theorem}

Choo and MacGillivray also showed that when $T$ is a tree, $k_{0}(T)=4$ if and
only if $T$ is a nontrivial odd star, and $k_{0}(T)=3$ otherwise. They also
showed that $k_{0}(C_{n})=4$ for each $n\geq3$. Celaya, Choo, MacGillivray,
and Seyffarth \cite{26CM-CCMS} continued the work of \cite{26CM-ChooMacG} and
considered complete bipartite graphs $K_{\ell,r}$. Since $\mathcal{C}_{2}(G)$
is disconnected for bipartite graphs, $k_{0}(K_{\ell,r})\geq3$. They proved
that equality holds if and only if $\ell$ and $r$ are both odd and that
$\mathcal{C}_{k}(K_{\ell,r})$ is Hamiltonian when $k\geq4$. Bard
\cite{26CM-Bard} expanded the latter result to complete multipartite graphs.

\begin{theorem}
\emph{\cite{26CM-Bard}}\hspace{0.1in}Fix $a_{1},...,a_{t}\in\mathbb{N}$. If
$k\geq2t$, then $\mathcal{C}_{k}(K_{a_{1},...,a_{t}})$ is Hamiltonian.
\end{theorem}

Bard improved this result for special cases by showing that $\mathcal{C}%
_{4}(K_{a_{1},a_{2},a_{3}})$ is Hamiltonian if and only if $a_{1}=a_{2}%
=a_{3}=1$, and, for $t\geq4$, $\mathcal{C}_{t+1}(K_{a_{1},...,a_{t}})$ is
Hamiltonian if and only if $a_{1}$ is odd and $a_{i}=1$ for $2\leq i\leq t$.
He showed that for each $k\geq4$ there exists a graph $G$ such that
$\mathcal{C}_{k}(G)$ is connected but not $2$-connected.

\begin{question}
\emph{\cite{26CM-Bard}}\hspace{0.1in}$(i)\hspace{0.1in}$Is $K_{2,2,2}$ the
only complete $3$-partite graph whose $5$-colouring graph is non-Hamiltonian?

$(ii)\hspace{0.1in}$Does there exist a connected $3$-colouring graph that is
not $2$-connected?

$(iii)\hspace{0.1in}$If $\mathcal{C}_{k}(G)$ is Hamiltonian, is $\mathcal{C}%
_{k+1}(G)$ always Hamiltonian?
\end{question}


Beier et al. \cite{26CM-Beier} considered the problem of determining which
graphs are realisable as colouring graphs. That is, given a graph $H$, when
does there exist a graph $G$ and an integer $k$ such that $H\cong%
\mathcal{C}_{k}(G)$? To this effect they determined that

\begin{itemize}
\item if $\mathcal{C}_{k}(G)$ is a complete graph, then it is $K_{k}$, and if
$k>1$ then $\mathcal{C}_{k}(G)=K_{k}$ if and only if $G=K_{1}$;

\item $K_{1}$ and $P_{2}$ are the only trees that are colouring graphs;

\item $C_{3},C_{4},C_{6}$ are the only cycles that are colouring graphs;

\item every tree is a subgraph of a colouring graph (thus there is no finite
forbidden subgraph characterisation of colouring graphs).
\end{itemize}


Other colouring graphs have also been considered. Haas \cite{26CM-Haas2012}
considered canonical and isomorphic colouring graphs. Two colourings of a
graph $G$ are \emph{isomorphic} if one results from permuting the names of the
colours of the other. A proper $k$-colouring of $G$ with colours $1,...,k$ is
\emph{canonical} with respect to an ordering $\pi=v_{1},...,v_{n}$ of the
vertices of $G$ if, for $1\leq c\leq k$, whenever colour $c$ is assigned to a
vertex $v_{i}$, each colour less than $c$ has been assigned to a vertex
$v_{j},\ j<i$. (Thus, a Grundy colouring $g$ becomes a canonical colouring if
we order the vertices of $G$ so that $v_{i}\prec v_{j}$ whenever
$g(v_{i})<g(v_{j})$.) For an ordering $\pi$ of the vertices of $G$, the
\emph{set of canonical }$k$\emph{-colourings of }$G$\emph{ under }$\pi$ is the
set $S_{\mathrm{Can}}(G)$ of pairwise nonisomorphic proper $k$-colourings of
$G$ that are lexicographically least under $\pi$. (Given colourings $c_{1}$
and $c_{2}$ of $G$ and an ordering $v_{1},...,v_{n}$ of $V(G)$, we say that
$c_{1}$ \emph{is lexicographically less than }$c_{2}$ if $c_{1}(v_{j}%
)<c_{2}(v_{j})$ for some integer $j,\ 1\leq j\leq n$, and $c_{1}(v_{i}%
)=c_{2}(v_{i})$ whenever $i<j$.) The \emph{canonical }$k$\emph{-colouring
graph} \textrm{C}$\mathrm{an}_{k}^{\pi}(G)$ is the graph with vertex set
$S_{\mathrm{Can}}(G)$ in which two colourings are adjacent if they differ at
exactly one vertex. Considering only nonisomorphic colourings, Haas defined
the \emph{isomorphic }$k$\emph{-colouring graph} $\mathcal{I}_{k}(G)$ to have
an edge between two colourings $c$ and $d$ if some representative of $c$
differs at exactly one vertex from some representative of $d$. Haas showed
that if the connected graph $G$ is not a complete graph, then \textrm{C}%
$\mathrm{an}_{k}^{\pi}(G)$ can be disconnected depending on the ordering $\pi$
and the difference $k-\chi(G)$. 

\begin{theorem}
\label{ThmCan}\emph{\cite{26CM-Haas2012}}\hspace{0.1in}$(i)\hspace{0.1in}$For
any connected graph $G\neq K_{n}$ and any $k\geq\chi(G)+1$ there exists an
ordering $\pi$ of $V(G)$ such that $\mathrm{Can}_{k}^{\pi}(G)$ is disconnected.

\begin{enumerate}
\item[$(ii)$] For any tree $T$ of order $n\geq4$ and any $k\geq3$ there is an
ordering $\pi$ of $V(T)$ such that $\mathrm{Can}_{k}^{\pi}(T)$ is Hamiltonian.

\item[$(iii)$] For any cycle $C_{n}$ and any $k\geq4$ there is an ordering
$\pi$ of $V(C_{n})$ such that $\mathrm{Can}_{k}^{\pi}(C_{n})$ is connected.
Moreover, $\mathrm{Can}_{3}^{\pi}(C_{4})$ and $\mathrm{Can}_{3}^{\pi}(C_{5})$
are connected for some $\pi$ but for all $n\geq6$, $\mathrm{Can}_{3}^{\pi
}(C_{n})$ is disconnected for all $\pi$.
\end{enumerate}
\end{theorem}

Haas and MacGillivray \cite{26CM-HM2018} extended this work and obtained a
variety of results on the connectedness and Hamiltonicity of the joins and
unions of graphs. They also obtained the following results.

\begin{theorem}
\emph{\cite{26CM-HM2018}\hspace{0.1in}}If $G$ is a bipartite graph on $n$
vertices, then there exists an ordering $\pi$ of $V(G)$ such that
$\mathrm{Can}_{k}^{\pi}(G)$ is connected for $k\geq n/2+1$.
\end{theorem}

\begin{theorem}
\emph{\cite{26CM-HM2018}\hspace{0.1in}}Let $G=K_{a_{1},...,a_{t}}$.

\begin{enumerate}
\item[$(i)$] For any $k\geq t$ there exists an ordering $\pi$ of $V(G)$ such
that $\mathrm{Can}_{k}^{\pi}(G)$ is connected.

\item[$(ii)$] If $a_{i}\geq2$ for each $i$, then for all vertex orderings
$\pi$ and $k\geq t+1$, $\mathrm{Can}_{k}^{\pi}(G)$ has a cut vertex and thus
is non-Hamiltonian, and if $t\geq3$, then $\mathrm{Can}_{k}^{\pi}(G)$ has no
Hamiltonian path.

\item[$(iii)$] For $t=2$, $K_{a_{1},a_{2}}$ has a vertex ordering $\pi$ such
that $\mathrm{Can}_{k}^{\pi}(K_{a_{1},a_{2}})$ has a Hamiltonian path for
$a_{1},a_{2}\geq2$ and $k\geq3$.
\end{enumerate}
\end{theorem}

Thus we see that all bipartite and complete multipartite graphs admit a vertex
ordering $\pi$ such that $\mathrm{Can}_{k}^{\pi}(G)$ is connected for large
enough values of $k$. Haas and MacGillivray also provided a vertex ordering
such that $\mathrm{Can}_{k}^{\pi}(G)$ is disconnected for all large values of
$k$.

Finbow and MacGillivray \cite{26CM-FM} studied the $k$-Bell colour graph and
the $k$-Stirling colour graph. The $k$-\emph{Bell colour graph} $\mathcal{B}%
_{k}(G)$ of $G$ is the graph whose vertices are the partitions of the vertices
of $G$ into at most $k$ independent sets, with different partitions $p_{1}$
and $p_{2}$ being adjacent if there is a vertex $x$ such that the restrictions
of $p_{1}$ and $p_{2}$ to $V(G)-\{x\}$ are the same partition. The
$k$-\emph{Stirling colour graph} $\mathcal{S}_{k}(G)$ of $G$ is the graph
whose vertices are the partitions of the vertices of $G$ into exactly $k$
independent sets, with adjacency as defined for $\mathcal{B}_{k}(G)$. They
showed, for example, that $\mathcal{B}_{n}(G)$ is Hamiltonian whenever $G$ is
a graph of order $n$ other than $K_{n}$ or $K_{n}-e$. As a consequence of
Theorem \ref{ThmCan}$(ii)$, $\mathcal{B}_{k}(T)$ is Hamiltonian whenever
$k\geq3$ and $T$ is a tree of order at least $4$, while $\mathcal{S}_{3}(T)$
has a Hamiltonian path. In addition, if $\mathcal{C}_{k}(G)$ is connected,
then so is $\mathcal{B}_{k}(G)$. They extended the result for $\mathcal{S}%
_{3}(T)$ to show that $\mathcal{S}_{k}(T)$ is Hamiltonian for any tree $T$ of
order $n\geq k+1$ and $k\geq4$.

Other variants of vertex colourings for which reconfiguration has been studied
include circular colourings \cite{26CM-BMMN, 26CM-BN}, acyclic colourings
\cite{26CM-Vai} and equitable colourings \cite{26CM-Vai}. Circular colourings
and $k$-colourings are special cases of homomorphisms, which we discuss in the
next subsection.

\subsection{Reconfiguration of Homomorphisms}

\label{SecWrochna}For graphs $G$ and $H$, a \emph{homomorphism} from $G$ to
$H$ is a mapping $\varphi:V(G)\rightarrow V(H)$ such that $\varphi
(u)\varphi(v)\in E(H)$ whenever $uv\in E(G)$. The collection of homomorphisms
from $G$ to $H$ is denoted by $\operatorname{Hom}(G,H)$. A $k$-colouring of
$G$ can be viewed as a homomorphism from $G$ to $K_{k}$. Thus we also refer to
a homomorphism from $G$ to $H$ as an $H$\emph{-colouring} of $G$. The
$H$\emph{-colouring graph} $\mathcal{C}_{H}(G)$ of $G$ has vertex set
$\operatorname{Hom}(G,H)$, and two homomorphisms are adjacent if one can be
obtained from the other by changing the colour of one vertex of $G$. For
$\alpha,\beta\in\operatorname{Hom}(G,H)$, an $\alpha,\beta$-walk in
$\mathcal{C}_{H}(G)$ is called an $H$-\emph{recolouring sequence} from
$\alpha$ to $\beta$. For a fixed graph $H$, the $H$-\emph{recolouring problem}
$H$\textsc{-Recolouring} is the problem of determining whether, given
$\alpha,\beta\in\operatorname{Hom}(G,H)$, there exists an $H$-recolouring
sequence from $\alpha$ to $\beta$. In the problem \textsc{Shortest}
$H$\textsc{-Recolouring}, one is also given an integer $\ell$, and the
question is whether the transformation can be done in at most $\ell$ steps.
Wrochna \cite{26CM-WrochnaDiss} approached the computational complexity of the
$H$-recolouring problem by using techniques from topology.

A graph $H$ has the \emph{monochromatic neighbourhood property (MNP)}, or is
an \emph{MNP-graph}, if for all pairs $a,b\in V(H)$, $|N_{H}(a)\cap
N_{H}(b)|\leq1$. Depending on whether $H$ has loops or not, MNP-graphs do not
contain $C_{4}$, or $K_{3}$ with one loop, or $K_{2}$ with both loops; $K_{3}$
and graphs with girth at least $5$ are all $C_{4}$-free. Note that
$3$-colourable graphs are MNP-graphs.

\begin{theorem}
\emph{\cite{26CM-WrochnaDiss}}\hspace{0.1in}If $H$ is an MNP-graph (possibly
with loops), then $H$\textsc{-Recolouring} and \textsc{Shortest}
$H$\textsc{-Recolouring} are in $P$.
\end{theorem}

Given positive integers $k$ and $q$ with $k\geq2q$, the \emph{circular clique}
$G_{k,q}$ has vertex set $\{0,1,...,k-1\}$, with $ij$ an edge whenever
$q\leq|i-j|\leq k-q$. A homomorphism $\varphi\in\operatorname{Hom}(G,G_{k,q})$
is called a \emph{circular colouring}. The \emph{circular chromatic number} of
$G$ is $\chi_{c}(G)=\inf\{k/q:\operatorname{Hom}(G,G_{k,q})\neq\emptyset\}$.
Brewster, McGuinness, Moore, and Noel \cite{26CM-BMMN} considered the
complexity of the $G_{k,q}$-recolouring problem.

\begin{theorem}
\emph{\cite{26CM-BMMN}}\hspace{0.1in}If $k$ and $q$ are fixed positive
integers with $k\geq2q$, then $G_{k,q}$\textsc{-Recolouring} is solvable in
polynomial time when $2\leq k/q<4$ and is PSPACE-complete for $k/q\geq4$.
\end{theorem}

The \emph{circular mixing number}\footnote{For comparison with $m_{0}(G)$ we
deviate slightly from the definition in \cite{26CM-BN} and adjust the results
accordingly.} of $G$, written $m_{c}(G)$, is $\inf\{r\in\mathbb{Q}:r\geq
\chi_{c}(G)$ and $\mathcal{C}_{G_{k,q}}(G)$ is connected whenever $k/q\geq
r\}$. Brewster and Noel \cite{26CM-BN} obtained bounds for $m_{c}(G)$ and
posed some interesting questions. They characterised graphs $G$ such that
$\mathcal{C}_{G}(G)$ is connected; this result requires a number of
definitions and we omit it here.

\begin{theorem}
\emph{\cite{26CM-BN}}\hspace{0.1in}$(i)\hspace{0.1in}$If $G$ is a graph of
order $n$, then $m_{c}(G)\leq2\operatorname{col}(G)$ and $m_{c}(G)\leq
\max\left\{  \frac{n+1}{2},m_{0}(G)\right\}  $. If $G$ has at least one edge,
then $m_{c}(G)\leq2\Delta(G)$.

\begin{enumerate}
\item[$(ii)$] If $G$ is a tree or a complete bipartite graph and $n\geq2$,
then $m_{c}(G)=2$.

\item[$(iii)$] If $G$ is nonbipartite, then $m_{c}(G)\geq\max\{4,\omega
(G)+1\}$.
\end{enumerate}
\end{theorem}

\begin{question}
\emph{\cite{26CM-BN}} 

\begin{enumerate}
\item[$(i)$] Is $m_{c}(G)$ always rational? When is it an integer?

\item[$(ii)$] Does there exist a real number $r$ such that $m_{c}(G)\leq
rm_{0}(G)$ for every graph $G$? If so, what is the smallest such $r$?
\end{enumerate}
\end{question}

\subsection{The $k$-Edge-Colouring Graph}

\label{SecKempe}In an attempt to prove the Four Colour Theorem, Alfred Bray
Kempe introduced the notion of changing map colourings by switching the
colours of regions in a maximal connected section of a map formed by regions
coloured with two specific colours, so as to eliminate a colour from regions
adjacent to an uncoloured region. (See e.g. \cite[Chapter 16]{26CM-CL}.) If we
consider proper edge-colourings of a graph $G$, then the subgraph $H$ of $G$
induced by all edges of two fixed colours has maximum degree 2; hence it
consists of the disjoint union of nontrivial paths and even cycles with edges
of alternating colours. These components of $H$ are now called
\emph{edge-Kempe chains}. We say that the proper $k$-edge-colourings $c_{1}$
and $c_{2}$ of $G$ are \emph{adjacent} in the $k$-edge-colouring graph
$\mathcal{EC}_{k}(G)$ if one can be obtained from the other by switching two
colours along an edge-Kempe chain. If a proper $k$-edge-colouring $c_{r}$ can
be converted to $c_{s}$ by a (possibly empty) sequence of edge-Kempe switches,
that is, if $c_{r}$ and $c_{s}$ are in the same component of $\mathcal{EC}%
_{k}(G)$, then we say that $c_{r}$ and $c_{s}$ are \emph{edge-Kempe
equivalent} and write $c_{r}\sim c_{s}$. Note that $\sim$ is an equivalence
relation; we may consider its equivalence classes on the set of $k$%
-edge-colourings of $G$. Two edge-colourings that differ only by a permutation
of colours are edge-Kempe equivalent, because the symmetric group $S_{k}$ is
generated by transpositions.

Most of the work on edge-Kempe equivalent edge-colourings has focused on the
number of equivalence classes of $k$-edge-colourings, i.e., the number of
components of $\mathcal{EC}_{k}(G)$, which we denote by $K^{\prime}(G,k)$. In
particular, the question of when $K^{\prime}(G,k)=1$ has received considerable
attention. In this section we allow our graphs to have multiple edges. We
denote the chromatic index (edge-chromatic number) of $G$ by $\chi^{\prime
}(G)$. Vizing (see e.g. \cite[Theorem 17.2]{26CM-CL}) proved that
$\Delta(G)\leq\chi^{\prime}(G)\leq\Delta(G)+1$ for any graph $G$.

Mohar \cite{26CM-Mohar} showed that if $k\geq\chi^{\prime}(G)+2$, then
$\mathcal{EC}_{k}(G)$ is connected, i.e., $K^{\prime}(G,k)=1$ for any graph
$G$, while if $G$ is bipartite and $k\geq\Delta(G)+1$, then $K^{\prime
}(G,k)=1$. He stated the characterisation of cubic bipartite graphs $G$ with
$K^{\prime}(G,3)=1$ as an open problem, and he conjectured that $K^{\prime
}(G,4)=1$ when $\Delta(G)\leq3$. (By K\"{o}nig's Theorem (see e.g.
\cite[Theorem 17.7]{26CM-CL}), $\chi^{\prime}(G)=3$ for a cubic bipartite
graph $G$.) McDonald, Mohar, and Scheide \cite{26CM-MMS} proved Mohar's
conjecture and showed that $K^{\prime}(K_{5},5)=6$.

\begin{theorem}
\emph{\cite{26CM-MMS}}\hspace{0.1in}$(i)\hspace{0.1in}$If $\Delta(G)\leq3$,
then $K^{\prime}(G,\Delta(G)+1)=1$.

\begin{enumerate}
\item[$(ii)$] If $\Delta(G)\leq4$, then $K^{\prime}(G,\Delta(G)+2)=1$.
\end{enumerate}
\end{theorem}

In \cite{26CM-belcasHaas2014}, belcastro and Haas provided partial answers to
Mohar's question on cubic bipartite graphs $G$ with $K^{\prime}(G,3)=1$. They
showed that all $3$-edge-colourings of planar bipartite cubic graphs are
edge-Kempe equivalent, and constructed infinite families of simple nonplanar
$3$-connected bipartite cubic graphs, all of whose $3$-edge-colourings are
edge-Kempe equivalent. In \cite{26CM-belcasHaas2017}, they investigated
$\mathcal{EC}_{k}(G)$ for $k$-edge-colourable $k$-regular graphs, and showed
that if such a graph is uniquely $k$-edge-colourable, then $\mathcal{EC}%
_{k}(G)$ is isomorphic to the Cayley graph of the symmetric group $S_{k}$ with
the set of all transpositions as generators.

\section{Reconfiguration of Dominating Sets}

\label{SecDom} There are several types of reconfiguration graphs of dominating
sets of a graph. Here we consider $k$-dominating graphs, $k$-total-dominating
graphs, and $\gamma$-graphs. In the first two cases, the vertices of the
reconfiguration graph correspond to (not necessarily minimal) dominating sets
of cardinality $k$ or less, whereas the vertices of $\gamma$-graphs correspond
to minimum dominating sets, also referred to as $\gamma$\emph{-sets}. A
minimal dominating set of maximum cardinality $\Gamma$ is called a $\Gamma
$\emph{-set}.

A graph $G$ is \emph{well-covered} if all its maximal independent sets have
cardinality $\alpha(G)$. A set $X\subseteq V(G)$ is \emph{irredundant} if each
vertex in $X$ dominates a vertex of $G$ (perhaps itself) that is not dominated
by any other vertex in $X$. An irredundant set is \emph{maximal irredundant}
if it has no irredundant proper superset. The lower and upper irredundant
numbers $\operatorname{ir}(G)$ and $\operatorname{IR}(G)$ of $G$ are,
respectively, the smallest and largest cardinalities of a maximal irredundant
set of $G$. If $X$ is a maximal irredundant set of cardinality
$\operatorname{ir}(G)$, we call $X$ an $\operatorname{ir}$-\emph{set}; an
$\operatorname{IR}$\emph{-set} is defined similarly.

A graph $G$ is \emph{irredundant perfect} if $\alpha(H)=\operatorname{IR}(H)$
for all induced subgraphs $H$ of $G$. Given a positive integer $k$, the family
$\mathcal{L}_{k}$ consists of all graphs $G$ containing vertices $x_{1}%
,\dots,x_{k}$ such that for each $i$, the subgraph induced by $N[x_{i}]$ is
complete, and $\{N[x_{i}]:1\leq i\leq k\}$ partitions $V(G)$. Let
$\mathcal{L}=\bigcup_{k\geq1}\mathcal{L}_{k}$. We use the graphs defined here
in the next section.

\subsection{The $k$-Dominating Graph}

\label{Sec_k_Dom}

The concept of $k$-dominating graphs was introduced by Haas and Seyffarth
\cite{26CM-HS2014} in 2014. This paper stimulated the work of Alikhani,
Fatehi, and Klav\v{z}ar \cite{26CM-davood}, Mynhardt, Roux, and Teshima
\cite{26CM-MRT}, Suzuki, Mouawad, and Nishimura \cite{26CM-SMN}, and their own
follow-up paper~\cite{26CM-HS2017}.

As is the case for $k$-colouring graphs, we seek to determine conditions for
the $k$-dominating graph $\mathcal{D}_{k}(G)$ to be connected. Haas and
Seyffarth \cite{26CM-HS2014} showed that any $\Gamma$-set $S$ of $G$ is an
isolated vertex of $\mathcal{D}_{\Gamma}(G)$ (because no proper subset of $S$
is dominating). Therefore, $\mathcal{D}_{\Gamma}(G)$ is disconnected whenever
$G$ has at least one edge (and thus at least two minimal dominating sets). In
particular, $\mathcal{D}_{n-1}(K_{1,n-1})$ is disconnected, but $\mathcal{D}%
_{k}(K_{1,n-1})$ is connected for all $k\in\{1,...,n\}-\{n-1\}$. This example
demonstrates that $\mathcal{D}_{k}(G)$ being connected does not imply that
$\mathcal{D}_{k+1}(G)$ is connected. However, Haas and Seyffarth showed that
if $k>\Gamma(G)$ and $\mathcal{D}_{k}(G)$ is connected, then $\mathcal{D}%
_{k+1}(G)$ is connected. They defined $d_{0}(G)$ to be the smallest integer
$\ell$ such that $\mathcal{D}_{k}(G)$ is connected for all $k\geq\ell$, and
noted that, for all graphs $G$, $d_{0}(G)$ exists because $\mathcal{D}_{n}(G)$
is connected. They bounded $d_{0}(G)$ as follows.

\begin{theorem}
\emph{\cite{26CM-HS2014}}\label{Thm_d0_bound}\hspace{0.1in}For any graph $G$
with at least one edge, $d_{0}(G)\geq\Gamma(G)+1$. If $G$ has at least two
disjoint edges, then $d_{0}(G)\leq\min\{n-1,\Gamma(G)+\gamma(G)\}$.
\end{theorem}

Haas and Seyffarth \cite{26CM-HS2017} showed that all independent dominating
sets of $G$ are in the same component of $D_{\Gamma(G)+1}(G)$ and established
the following upper bound for $d_{0}(G)$; for a graph with $\gamma=\alpha$ it
improves the bound in Theorem \ref{Thm_d0_bound}.

\begin{theorem}
\emph{\cite{26CM-HS2017}}\label{Thm_d0_bound2}\hspace{0.1in}For any graph $G$,
$d_{0}(G)\leq\Gamma(G)+\alpha(G)-1$. Furthermore, if $G$ is triangle-free,
then $d_{0}(G)\leq\Gamma(G)+\alpha(G)-2$.
\end{theorem}

Graphs for which equality holds in the lower bound in Theorem
\ref{Thm_d0_bound} (provided they are connected and nontrivial) include
bipartite graphs, chordal graphs \cite{26CM-HS2014}, graphs with $\alpha\leq
2$, graphs that are perfect and irredundant perfect, well-covered graphs with
neither $C_{4}$ nor $C_{5}$ as subgraph, well-covered graphs with girth at
least five, well-covered claw-free graphs without $4$-cycles, well-covered
plane triangulations, and graphs in the class $\mathcal{L}$ \cite{26CM-HS2017}.

\begin{figure}[h]
\begin{center}
\scalebox{0.8}{\includegraphics{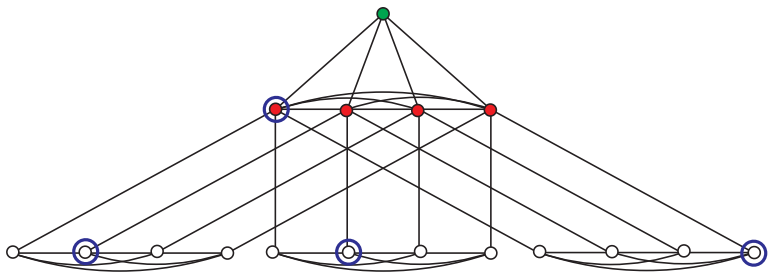}}
\end{center}
\caption[Short figure caption]{A graph $G$ with $\gamma(G)=\Gamma(G)=4$ and
$d_{0}(G)=7=\Gamma(G)+\gamma(G)-1$}%
\label{Figk4r3}%
\end{figure}

\begin{figure}[h]
\begin{center}
\scalebox{0.8}{\includegraphics{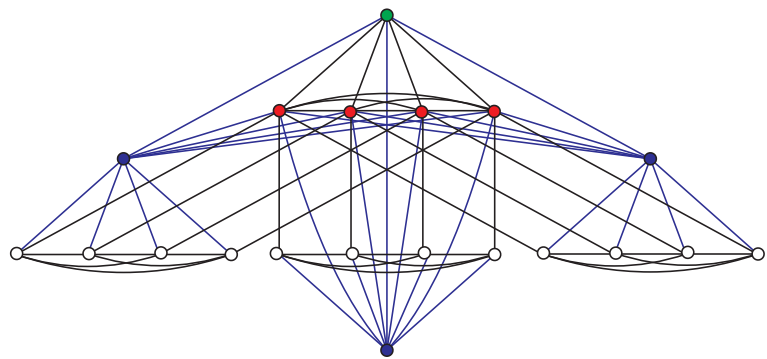}}
\end{center}
\caption[Short figure caption]{A graph $Q$ with $\gamma(Q)=3$, $\Gamma(Q)=4$
and $d_{0}(Q)=7=\Gamma(Q)+\gamma(Q)$}%
\label{Figk4+3}%
\end{figure}

Suzuki et al.~\cite{26CM-SMN} were first to exhibit graphs for which
$d_{0}>\Gamma+1$. They constructed an infinite class of graphs $G_{(d,b)}$ (of
tree-width $2b-1$) for which $d_{0}(G_{(d,b)})=\Gamma(G_{(d,b)})+2$; the
smallest of these is $G_{(2,3)}\cong P_{3}\boksie K_{3}$, which is planar.
Haas and Seyffarth \cite{26CM-HS2017} also found a graph $G_{4}$ such that
$d_{0}(G_{4})=\Gamma(G_{4})+2$, and they mentioned that they did not know of
the existence of any graphs with $d_{0}>\Gamma+2$. Mynhardt et
al.~\cite{26CM-MRT} constructed classes of graphs that demonstrate (a) the
existence of graphs with arbitrary upper domination number $\Gamma\geq3$,
arbitrary domination number in the range $2\leq\gamma\leq\Gamma$, and
$d_{0}=\Gamma+\gamma-1$ (see Figure \ref{Figk4r3} for an example), and (b) the
existence of graphs with arbitrary upper domination number $\Gamma\geq3$,
arbitrary domination number in the range $1\leq\gamma\leq\Gamma-1$, and
$d_{0}=\Gamma+\gamma$ (see Figure \ref{Figk4+3} for an example). For
$\gamma\geq2$, this was the first construction of graphs with $d_{0}%
=\Gamma+\gamma$. These results are best possible in both cases, since it
follows from Theorems \ref{Thm_d0_bound} and \ref{Thm_d0_bound2} that
$d_{0}(G)\leq\min\{\Gamma(G)+\gamma(G),2\Gamma(G)-1\}$ for any graph $G$.

Suzuki et al.~\cite{26CM-SMN} related the connectedness of $\mathcal{D}%
_{k}(G)$ to matchings in $G$ by showing that if $G$ has a matching of size (at
least) $\mu+1$, then $\mathcal{D}_{n-\mu}(G)$ is connected. This result is
best possible with respect to the size of a maximum matching, since the path
$P_{2k}$ has matching number $\mu=k=\Gamma(P_{2k})=n-\mu$, hence
$\mathcal{D}_{n-\mu}(P_{2k})$ is disconnected. It also follows that the
diameter of $\mathcal{D}_{n-\mu}(G)$ is in $O(n)$ for a graph $G$ with a
matching of size $\mu+1$. On the other hand, they constructed an infinite
family of graphs $G_{n}$ of order $63n-6$ such that $\mathcal{D}_{\gamma
(G)+1}(G_{n})$ has exponential diameter $\Omega(2^{n})$.

\begin{question}
$\null$

\begin{enumerate}
\item[$(i)$] \emph{\cite{26CM-HS2014}}\hspace{0.1in}Characterise graphs for
which $d_{0}=\Gamma+1$.

\item[$(ii)$] \emph{\cite{26CM-MRT}}\hspace{0.1in}Is it true that
$d_{0}(G)=\Gamma(G)+1$ when $G$ is triangle-free?

\item[$(iii)$] \emph{\cite{26CM-HS2014}}\hspace{0.1in}When is $\mathcal{D}%
_{k}(G)$ Hamiltonian?

\item[$(iv)$] \emph{\cite{26CM-MRT}}\hspace{0.1in}Suppose $\mathcal{D}_{i}(G)$
and $\mathcal{D}_{j}(G)$ are connected and $i<j$. How are $\operatorname{diam}%
(\mathcal{D}_{i}(G))$ and $\operatorname{diam}(\mathcal{D}_{j}(G))$ related?
(If $i>\Gamma(G)$, then $\operatorname{diam}(\mathcal{D}_{i}(G))\geq
\operatorname{diam}(\mathcal{D}_{j}(G))$.)
\end{enumerate}
\end{question}


Haas and Seyffarth~\cite{26CM-HS2014} considered the question of which graphs
are realisable as $k$-dominating graphs and observed that for $n\geq4$,
$\mathcal{D}_{2}(K_{1,n-1})=K_{1,n-1}$. Alikhani et al.~\cite{26CM-davood}
proved that these stars are the only graphs with this property, i.e. if $G$ is
a graph of order $n$ with no isolated vertices such that $n\geq2$, $\delta
\geq1$, and $G\cong\mathcal{D}_{k}(G)$, then $k=2$ and $G\cong K_{1,n-1}$ for
some $n\geq4$. They also showed that $C_{6},C_{8},P_{1}$ and $P_{3}$ are the
only cycles or paths that are dominating graphs of connected graphs
($\mathcal{D}_{2}(K_{3})=C_{6}$, $\mathcal{D}_{3}(P_{4})=C_{8}$,
$\mathcal{D}_{1}(K_{1})=P_{1}$ and $\mathcal{D}_{2}(K_{2})=P_{3}$). They
remarked that $\mathcal{D}_{n}(G)$ has odd order for every graph $G$ (since
$G$ has an odd number of dominating sets \cite{26CM-BCS}), and showed that if
$m$ is odd and $0<m<2^{n}$, then there exists a graph $X$ of order $n$ such
that $\mathcal{D}_{n}(X)$ has order $m$.

It is obvious that $\mathcal{D}_{k}(G)$ is bipartite for any graph $G$ of
order $n$ and any $k$ such that $\gamma(G)\leq k\leq n$; in fact,
$\mathcal{D}_{k}(G)$ is an induced subgraph of $Q_{n}-v$, a hypercube with one
vertex deleted~\cite{26CM-davood}.

\begin{question}
Which induced subgraphs of $Q_{n}$ occur as $\mathcal{D}_{k}(G)$ for some
$n$-vertex graph $G$ and some integer $k$?
\end{question}

\subsection{The $k$-Total-Dominating Graph}

\label{Sec_k_TotalDom}

For a graph $G$ without isolated vertices, a set $S\subseteq V(G)$ is a
\emph{total-dominating set (TDS) }if every vertex of $G$ is adjacent to a
vertex in $S$. We denote the minimum (maximum, respectively) cardinality of a
minimal TDS by $\gamma_{t}(G)$ ($\Gamma_{t}(G)$, respectively). Alikhani,
Fatehi, and Mynhardt \cite{26CM-AFM} initiated the study of \emph{$k$%
}-total-dominating graphs (see Section~\ref{26CM-Sec_Intro}). Since any TDS is
a dominating set, $\mathcal{D}_{k}^{t}(G)$ is an induced subgraph of
$\mathcal{D}_{k}(G)$ for any isolate-free graph $G$ and any integer
$k\geq\gamma_{t}(G)$. However, since $\Gamma$ and $\Gamma_{t}$ are not
comparable (for $n$ large enough, $\Gamma_{t}(K_{1,n})=2<\Gamma(K_{1,n})=n$
but $\Gamma(P_{n})<\Gamma_{t}(P_{n})$), the two graphs $\mathcal{D}_{k}(G)$
and $\mathcal{D}_{k}^{t}(G)$ can be different.

To study the connectedness of $\mathcal{D}_{k}^{t}(G)$, we define $d_{0}%
^{t}(G)$ similar to $d_{0}(G)$ (Section~\ref{Sec_k_Dom}). Unlike
$\mathcal{D}_{\Gamma}(G)$, there are nontrivial connected graphs $G$ such that
$\mathcal{D}_{\Gamma_{t}}(G)$ is connected and $d_{0}^{t}(G)=\Gamma_{t}(G)$,
as shown below. The unique neighbour of a vertex of degree one is called a
\emph{stem}. Denote the set of stems of $G$ by $S(G)$.

\begin{theorem}
\label{Thm_TDS}\emph{\cite{26CM-AFM}}\hspace{0.1in}If $G$ is a connected graph
of order $n\geq3$, then

\begin{enumerate}
\item[$(i)$] $\mathcal{D}_{\Gamma_{t}}^{t}(G)$ is connected if and only if
$S(G)$ is a TDS of $G$,

\item[$(ii)$] $\Gamma_{t}(G)\leq d_{0}^{t}(G)\leq n$,

\item[$(iii)$] any isolate-free graph $H$ is an induced subgraph of a graph
$G$ such that $\mathcal{D}_{\Gamma_{t}}^{t}(G)$ is connected ($G$ is the
corona of $H$),

\item[$(iv)$] if $G$ is a connected graph of order $n\geq3$ such that $S(G)$
is a TDS, then $\mathcal{D}_{\gamma_{t}}^{t}(G)$ is connected ($S(G)$ is the
unique TDS).
\end{enumerate}
\end{theorem}

The lower bound in Theorem \ref{Thm_TDS}$(ii)$ is realised if and only if $G$
has exactly one minimal TDS, i.e. if and only if $S(G)$ is a TDS. The upper
bound is realised if and only if $\Gamma_{t}(G)=n-1$, i.e. if and only if $n$
is odd and $G$ is obtained from $\frac{n-1}{2}K_{2}$ by joining a new vertex
to at least one vertex of each $K_{2}$.

For specific graph classes, Alikhani et al.~\cite{26CM-AFM} showed that
$d_{0}^{t}(C_{n})=\Gamma_{t}(C_{n})+1$ if $n\neq8$, while if $n=8$, then
$d_{0}^{t}(C_{8})=\Gamma_{t}(C_{8})+2$. Hence $\mathcal{D}_{\Gamma_{t}+1}%
^{t}(C_{8})$ is disconnected, making $C_{8}$ the only known graph with this
property. For paths, $d_{0}^{t}(P_{2})=\Gamma_{t}(P_{2})=d_{0}^{t}%
(P_{4})=\Gamma_{t}(P_{4})=2$ and $d_{0}^{t}(P_{n})=\Gamma_{t}(P_{n})+1$ if
$n=3$ or $n\geq5$.

As shown in \cite{26CM-AFM}, $Q_{n}$ and $K_{1,n}$, $n\geq2$, are realisable
as total-dominating graphs, and $C_{4},C_{6},C_{8},C_{10},P_{1},P_{3}$ are the
only realisable cycles and paths.

\begin{question}
\emph{\cite{26CM-AFM}}

\begin{enumerate}
\item[$(i)$] Construct classes of graphs $G_{r}$ such that $d_{0}^{t}%
(G_{r})-\Gamma_{t}(G_{r})\geq r\geq2$.

\item[$(ii)$] Find more classes of graphs that can/cannot be realised as
$k$-total-domination graphs.

\item[$(iii)$] Note that $\mathcal{D}_{3}^{t}(P_{3})\cong P_{3}$. Characterise
graphs $G$ such that $\mathcal{D}_{k}^{t}(G)\cong G$ for some~$k$.
\end{enumerate}
\end{question}

\subsection{Jump $\gamma$-Graphs}

\label{SecJump}Sridharan and Subramanaian \cite{26CM-Sub} introduced jump
$\gamma$-graphs $\mathcal{J}(G,\gamma)$ in 2008; they used the notation
$\gamma\cdot G$ instead of $\mathcal{J}(G,\gamma)$. The $\gamma$-graphs
$\mathcal{J}(G,\gamma)$ for $G\in\{P_{n},C_{n}\}$ were determined in
\cite{26CM-Sub}, as were the graphs $\mathcal{J}(H_{k,n},\gamma)$ for some
values of $k$ and $n$, where $H_{k,n}$ is a \emph{Harary graph}, i.e. a
$k$-connected graph of order $n$ and minimum possible size $\left\lceil
kn/2\right\rceil $. The authors of \cite{26CM-Sub} showed that if $T$ is a
tree, then $\mathcal{J}(T,\gamma)$ is connected. Haas and Seyffarth
\cite{26CM-HS2014} showed that if $\mathcal{D}_{\gamma(G)+1}(G)$ is connected,
then $\mathcal{J}(G,\gamma)$ is connected, thus relating $k$-dominating graphs
to $\gamma$-graphs.

Sridharan and Subramanaian \cite{26CM-Srid} showed that trees and unicyclic
graphs can be realised as jump $\gamma$-graphs. Denoting the graph obtained by
joining the two vertices of $K_{2,3}$ of degree $3$ by $\Delta_{3}$, they
showed that if $H$ contains $\Delta_{3}$ as an induced subgraph, then $H$ is
not realisable as a $\gamma$-graph $\mathcal{J}(G,\gamma)$. Following the same
line of enquiry, Lakshmanan and Vijayakumar \cite{26CM-LV} proved that if $H$
is a $\gamma$-graph, then $H$ contains none of $K_{2,3},K_{2}\vee P_{3}%
,(K_{1}\cup K_{2})\vee2K_{1}$ as an induced subgraph. They showed that the
collection of $\gamma$-graphs is closed under the Cartesian product and that a
disconnected graph is realisable if and only if all its components are
realisable. They also proved that if $G$ is a connected cograph, then
$\operatorname{diam}(\mathcal{J}(G,\gamma))\leq2$, where $\operatorname{diam}%
(\mathcal{J}(G,\gamma))=1$ if and only if $G$ has a universal vertex. Bie\'{n}
\cite{26CM-Bien} studied $\mathcal{J}(T,\gamma)$ for trees of diameter at most
$5$ and for certain caterpillars.

In his Master's thesis \cite{26CM-Dyck}, Dyck illustrated a connection between
$\gamma$-graphs and Johnson graphs. The \emph{Johnson graph} $J(n,k)$ is the
graph whose vertex set consists of all $k$-subsets of $\{1,...,n\}$, where two
vertices are adjacent whenever their corresponding sets intersect in exactly
$k-1$ elements.

\begin{theorem}
\emph{\cite{26CM-Dyck}}\hspace{0.1in}A graph $H$ is realisable as
$\mathcal{J}(G,\gamma)$, where $G$ is an $n$-vertex graph with $\gamma(G)=k$,
if and only if $H$ is isomorphic to an induced subgraph of $J(n,k)$.
\end{theorem}

Edwards, MacGillivray, and Nasserasr \cite{26CM-EMN} obtained results which
hold for jump and slide $\gamma$-graphs; we report their results in Theorem
\ref{Theorem_EMN}.

\subsection{Slide $\gamma$-Graphs}

\label{SecSlide}

Fricke, Hedetniemi, Hedetniemi, and Hutson \cite{26CM-FHHH} introduced slide
$\gamma$-graphs $\mathcal{S}(G,\gamma)$ in 2011; they used the notation
$G(\gamma)$ instead of $\mathcal{S}(G,\gamma)$. They showed that every tree is
realisable as a slide $\gamma$-graph, that $\mathcal{S}(T,\gamma)$ is
connected and bipartite if $T$ is a tree, and that $\mathcal{S}(G,\gamma)$ is
triangle-free if $G$ is triangle-free. They determined $\mathcal{S}(G,\gamma)$
for a number of graph classes, including complete and complete bipartite
graphs, paths and cycles.

Connelly, Hedetniemi, and Hutson \cite{26CM-ConHH} extended the realisability
result obtained in \cite{26CM-FHHH}.

\begin{theorem}
\emph{\cite{26CM-ConHH}}\hspace{0.1in}Every graph is realisable as a $\gamma
$-graph $\mathcal{S}(G,\gamma)$ of infinitely many graphs~$G$.
\end{theorem}

Connelly et al. \cite{26CM-ConHH} also showed that the $\gamma$-graphs of all
graphs of order at most $5$ are connected and characterised graphs of order
$6$ with disconnected $\gamma$-graphs.

Edwards et al.~\cite{26CM-EMN} investigated the order, diameter, and maximum
degree of jump and slide $\gamma$-graphs of trees, providing answers to
questions posed in \cite{26CM-FHHH}.

\begin{theorem}
\label{Theorem_EMN}\emph{\cite{26CM-EMN}}\hspace{0.1in}If $T$ is a tree of
order $n$ having $s$ stems, then

\begin{enumerate}
\item[$(i)$] $\Delta(\mathcal{S}(T,\gamma))\leq n-\gamma(T)$ and
$\Delta(\mathcal{J}(T,\gamma))\leq n-\gamma(T)$,

\item[$(ii)$] $\operatorname{diam}(\mathcal{S}(T,\gamma))\leq2(2\gamma(T)-s)$
and $\operatorname{diam}(\mathcal{J}(T,\gamma))\leq2\gamma(T)$,

\item[$(iii)$] $|V(\mathcal{S}(T,\gamma))|=|V(\mathcal{J}(T,\gamma
))|\leq((1+\sqrt{13})/2)^{\gamma(T)}$.
\end{enumerate}
\end{theorem}

It follows that the maximum degree and diameter of $\gamma$-graphs of trees
are linear in $n$. Edwards et al. exhibited an infinite family of trees to
demonstrate that the bounds in Theorem \ref{Theorem_EMN}$(i)$ are sharp and
mentioned that there are no known trees for which $\operatorname{diam}%
(\mathcal{S}(T,\gamma))$ or $\operatorname{diam}(\mathcal{J}(T,\gamma))$
exceeds half the bound given in Theorem \ref{Theorem_EMN}$(ii)$. They also
demonstrated that $|V(\mathcal{S}(T,\gamma))|>2^{\gamma(T)}$ for infinitely
many trees.

\begin{question}
$\null$

\begin{enumerate}
\item[$(i)$] \emph{\cite{26CM-FHHH}}
Which graphs are $\gamma$-graphs of trees?

\item[$(ii)$] \emph{\cite{26CM-MT}}\hspace{0.1in}Is every bipartite graph the
$\gamma$-graph of a \textbf{bipartite} graph?
\end{enumerate}
\end{question}

\subsection{Irredundance}

Mynhardt and Teshima \cite{26CM-MT} studied slide reconfiguration graphs for
other domination parameters. In particular, for an arbitrary given graph $H$
they constructed a graph $G_{H}$ to show that $H$ is realisable as the slide
$\Gamma$-graph $\mathcal{S}(G_{H},\Gamma)$ of $G_{H}$. Although $G_{H}$
satisfies $\Gamma(G_{H})=\operatorname{IR}(G_{H})$, it has more
$\operatorname{IR}$-sets than $\Gamma$-sets. Hence $H$ is not an
$\operatorname{IR}$-graph of $G_{H}$. They left the problem of whether all
graphs are $\operatorname{IR}$-graphs open. Mynhardt and Roux \cite{26CM-MR}
responded as follows.

\begin{theorem}
\emph{\cite{26CM-MR}}\hspace{0.1in}$(i)\hspace{0.1in}$All disconnected graphs
can be realised as $\operatorname{IR}$-graphs.

\begin{enumerate}
\item[$(ii)$] Stars $K_{1,k}$ for $k\geq2$, the cycles $C_{5},C_{6},C_{7}$,
and the paths $P_{3},P_{4},P_{5}$ are not $\operatorname{IR}$-graphs.
\end{enumerate}
\end{theorem}

Mynhardt and Roux also showed that the double star $S(2,2)$ (obtained by
joining the central vertices of two copies of $P_{3}$), and the tree obtained
by joining a new leaf to a leaf of $S(2,2)$, are the unique smallest
$\operatorname{IR}$-trees with diameters $3$ and $4$, respectively. The only
connected $\operatorname{IR}$-graphs of order $4$ are $K_{4}$ and $C_{4}$. We
close with one of their questions and a conjecture.

\begin{conjecture}
\emph{\cite{26CM-MR}}\hspace{0.1in}$P_{n}$ is not an $\operatorname{IR}$-graph
for each $n\geq3$, and $C_{n}$ is not an $\operatorname{IR}$-graph for each
$n\geq5$.
\end{conjecture}

\begin{question}
\emph{\cite{26CM-MR}}\hspace{0.1in}Are complete graphs and $C_{4}$ the only
claw-free $\operatorname{IR}$-graphs?
\end{question}

\noindent\textbf{{\Large {\large Acknowledgemen}\textbf{{\Large {\large t}}}%
}\hspace{0.1in}}This survey was published as \cite{MN2019}.

\label{References}


\begin{thebibliography}{99}                                                                                               %


\bibitem {26CM-davood}S. Alikhani, D. Fatehi, and S. Klav\v{z}ar. On the
structure of dominating graphs. \emph{Graphs Combin.}, 33:665--672,\ 2017.

\bibitem {26CM-AFM}S. Alikhani, D. Fatehi, and C. M. Mynhardt. On $k$-total
dominating graphs. \emph{Australas.J.Combin.}, 73:313--333, 2019.

\bibitem {26CM-Bard}S. Bard. \emph{Gray code numbers of complete multipartite
graphs}. Master's thesis, University of Victoria, 2014. http://hdl.handle.net/1828/5815

\bibitem {26CM-Beier}J. Beier, J. Fierson, R. Haas, H. M. Russell, Heather,
and K. Shavo. Classifying coloring graphs.\emph{ Discrete Math.},
339:2100--2112, 2016.

\bibitem {26CM-belcasHaas2014}s. m. belcastro and R. Haas. Counting
edge-Kempe-equivalence classes for $3$-edge-colored cubic graphs.
\emph{Discrete Math}., 325:77--84, 2014.

\bibitem {26CM-belcasHaas2017}s. m. belcastro and R. Haas.
Edge-Kempe-equivalence graphs of class-$1$ regular graphs. \emph{Australas. J.
Combin}., 69:197--214, 2017.

\bibitem {26CM-Bien}A. Bie\'{n}. Gamma graphs of some special classes of
trees. \emph{Ann. Math. Sil.}, 29:25--34, 2015.

\bibitem {26CM-BB}M. Bonamy and N. Bousquet. Recoloring graphs via tree
decompositions. \emph{European J. Combin}., 69:200--213, 2018.

\bibitem {26CM-BJLPP}M. Bonamy, M. Johnson, I. Lignos, V. Patel, and D.
Paulusma. On the diameter of reconfiguration graphs for vertex colourings.
\emph{Electron. Notes Discrete Math.}, 38:161--166, 2011.

\bibitem {26CM-BC}P. Bonsma and L. Cereceda. Finding paths between graph
colourings: PSPACE-completeness and superpolynomial distances. \emph{Theor.
Comput. Sci.}, 410:5215--5226,\textbf{ }2009.

\bibitem {26CM-BCHJ}P. Bonsma, L. Cereceda, J. van den Heuvel, and M. Johnson.
Finding paths between graph colourings: computational complexity and possible
distances. \emph{Electron. Notes Discrete Math.}, 29:463--469, 2007.

\bibitem {26CM-BM}P. Bonsma and A. E. Mouawad. The complexity of bounded
length graph recoloring. Manuscript. Arxiv.org/pdf/1404.0337.pdf

\bibitem {26CM-BMNR}P. Bonsma, A. E. Mouawad, N. Nishimura, and V. Raman The
complexity of bounded length graph recoloring and CSP reconfiguration. In
\emph{Proceedings of the 9th International Symposium on Parameterized and
Exact Computation}, IPEC 2014, Wroclaw, Poland, pp. 110--121, 2014.

\bibitem {26CM-BCS}A. E. Brouwer, P. Csorba, and A. Schrijver: The number of
dominating sets of a finite graph is odd. Manuscript, 2009.

\bibitem {26CM-BMMN}R. C. Brewster, S. McGuinness, B. Moore, and J. Noel. A
dichotomy theorem for circular colouring reconfiguration.\emph{ Theoret.
Comput. Sci.}, 639:1--13, 2016.

\bibitem {26CM-BN}R. C. Brewster and J. A. Noel. Mixing Homomorphisms,
Recolorings, and Extending Circular Precolorings. \emph{J. Graph Theory},
80:173--198, 2015.

\bibitem {26CM-CCMS}M. Celaya, K. Choo, G. MacGillivray, and K. Seyffarth.
Reconfiguring $k$-colourings of complete bipartite graphs.\emph{ Kyungpook
Math. J.}, 56:647--655, 2016.

\bibitem {26CM-CHJ1}L. Cereceda, J. van den Heuvel, and M. Johnson.
Connectedness of the graph of vertex-colourings. \emph{Discrete Math.},
308:913--919, 2008.

\bibitem {26CM-CHJ2}L. Cereceda, J. van den Heuvel, and M. Johnson. Mixing
$3$-colourings in bipartite graphs. \emph{European J. Combin}., 30:1593--1606, 2009.

\bibitem {26CM-CHJ3}L. Cereceda, J. van den Heuvel, and M. Johnson. Finding
paths between $3$-colorings. \emph{J. Graph Theory}, 67:69--82, 2011.

\bibitem {26CM-CL}G. Chartrand, L. Lesniak, and P. Zhang. \emph{Graphs
\&\ Digraphs},\emph{ }6th ed. Chapman and Hall/CRC, Boca Raton, 2016.

\bibitem {26CM-ChooThesis}K. Choo. \emph{The existence of grey codes for
proper }$k$\emph{-colourings of graphs}. Master's thesis, University of
Victoria, 2002.

\bibitem {26CM-ChooMacG}K. Choo and G. MacGillivray. Gray code numbers for
graphs. \emph{Ars Math. Contemp.}, 4:125--139, 2011.

\bibitem {26CM-ConHH}E. Connelly, S.T. Hedetniemi, and K.R. Hutson. A note on
$\gamma$-Graphs. \emph{AKCE Intr. J. Graphs Comb}., 8:23--31, 2010.

\bibitem {26CM-Dyck}A. R. J. Dyck. \emph{The realisability of }$\gamma
$\emph{-graphs}. Master's thesis, Simon Fraser University, 2017. summit.sfu.ca/item/17513

\bibitem {26CM-DFFV}M. Dyer, A. D. Flaxman, A. M. Frieze, and E. Vigoda.
Randomly coloring sparse random graphs with fewer colors than the maximum
degree. \emph{Random Structures Algorithms}, 29:450--465,\textbf{ }2006.

\bibitem {26CM-DGM}M. Dyer, C. Greenhill, and M. Molloy. Very rapid mixing of
the Glauber dynamics for proper colorings on bounded-degree graphs.
\emph{Random Structures Algorithms}, 20:98--114, 2001.

\bibitem {26CM-EMN}M. Edwards, G. MacGillivray, and S.
Nasserasr.\ Reconfiguring minimum dominating sets: the $\gamma$-graph of a
tree. \emph{Discuss. Math. Graph Theory}, 38:703--716, 2018.

\bibitem {26CM-FJP}C. Feghali, M. Johnson, and D. Paulusma. A reconfigurations
analogue of Brooks' theorem and its consequences. \emph{J. Graph Theory},
83:340--358, 2016.

\bibitem {26CM-FM}S. Finbow and G. MacGillivray. Hamiltonicity of Bell and
Stirling colour graphs. Manuscript, 2014.

\bibitem {26CM-FHHH}G. H. Fricke, S. M. Hedemiemi, S. T. Hedetniemi, and K. R.
Hutson. $\gamma$-Graphs of graphs.\emph{ Discuss. Math. Graph Theory},
31:517--531, 2011.

\bibitem {26CM-Haas2012}R. Haas. The canonical coloring graph of trees and
cycles.\emph{ Ars Math. Contemp.}, 5:149--157, 2012.

\bibitem {26CM-HM2018}R. Haas and G. MacGillivray. Connectivity and
Hamiltonicity of canonical colouring graphs of bipartite and complete
multipartite graphs. \emph{Algorithms} (Basel, Paper No. 40), 11, 14 pp., 2018.

\bibitem {26CM-HS2014}R. Haas and K. Seyffarth. The $k$-dominating graph.
\emph{Graphs Combin.}, 30:609--617, 2014.

\bibitem {26CM-HS2017}R. Haas and K. Seyffarth. Reconfiguring dominating sets
in some well-covered and other classes of graphs. \emph{Discrete Math.}
340:1802--1817, 2017.

\bibitem {26CM-HIMNOST}A. Haddadan, T. Ito, A. E.Mouawad, N. Nishimura, H.
Ono, A. Suzuki, and Y. Tebbal. The complexity of dominating set
reconfiguration. \emph{Theoret. Comput. Sci.}, 651:37--49, 2016.

\bibitem {26CM-Ito2}T. Ito, E. D. Demaine, N. J. A. Harvey, C. H.
Papadimitriou, M. Sideri, R. Uehara, and Y. Uno. On the complexity of
reconfiguration problems. \emph{Theoret. Comput. Sci.} 412:1054--1065, 2011.

\bibitem {26CM-IKD}T. Ito, M. Kaminski, and E. D. Demaine. Reconfiguration of
list edge-colorings in a graph. In \emph{Algorithms and data str\emph{u}%
ctures, LNCS}, 5664, Springer, Berlin, pp. 375--386, 2009.

\bibitem {26CM-Jerrum}M. Jerrum. A very simple algorithm for estimating the
number of $k$-colorings of a low-degree graph. \emph{Random Structures
Algorithms}, 7:157--165, 1995.

\bibitem {26CM-JKKPP}M. Johnson, D. Kratsch, S. Kratsch, V. Patel, and D.
Paulusma. Finding shortest paths between graph colourings. \emph{Parameterized
and exact computation}. In \emph{Lecture Notes in Comput. Sci.} 8894,
Springer, Cham, pp. 221--233, 2014.

\bibitem {26CM-JKKPP2}M. Johnson, D. Kratsch, S. Kratsch, V. Patel, and D.
Paulusma. Finding shortest paths between graph colourings. \emph{Algorithmica}%
, 75:295--321, 2016.

\bibitem {26CM-LV}S. A. Lakshmanan and A. Vijayakumar. The gamma graph of a
graph. \emph{AKCE Intr. J. Graphs Comb}., 7:53--59,\textbf{ }2010.

\bibitem {26CM-LMPRS}D. Lokshtanov, A. E. Mouawad, F. Panolan, M. S.
Ramanujan, and S. Saurabh. Reconfiguration on Sparse Graphs. In
\emph{Proceedings of the 14th International Symposium on Algorithms and Data
Structures}, WADS 2015, Victoria, BC, Canada, pp. 506--517, 2015.

\bibitem {26CM-LM}B. Lucier and M. Molloy. The Glauber dynamics for colorings
of bounded degree trees. \emph{SIAM J. Discrete Math.}, 25:827--853, 2011.

\bibitem {26CM-MMS}J. McDonald, B. Mohar, and D. Scheide. Kempe equivalence of
edge-colorings in subcubic and subquartic graphs. \emph{J. Graph Theory,
}70:226--239, 2012.

\bibitem {26CM-Mohar}B. Mohar. Kempe equivalence of colorings. \emph{Graph
theory in Paris}, 287--297. Trends Math., Birkh\"{a}user, Basel, 2007.

\bibitem {26CM-Molloy}M. Molloy. The glauber dynamics on colorings of a graph
with high girth and maximum degree. \emph{SIAM J. Comput.},\emph{
}33:721--737, 2004.

\bibitem {26CM-MNRSS}A. E. Mouawad, N. Nishimura, V. Raman, N. Simjour, and A.
Suzuki. On the parameterized complexity of reconfiguration problems.
\emph{Algorithmica}, 78:274--297, 2017.

\bibitem {MN2019}C.~M.~Mynhardt and S.~Nasserasr, Reconfiguration of
colourings and dominating sets in graphs, in F. Chung, R. Graham, F.
Hoffman,L. Hogben, Ron Mullin and Doug West (Eds.), \emph{50 Years of
Combinatorics, Graph Theory, and Computing}, Chapman and Hall/CRC Press, 2019.

\bibitem {26CM-MR}C. M. Mynhardt and A. Roux. Irredundance graphs. Manuscript,
2018. arXiv:1812.03382v1

\bibitem {26CM-MRT}C. M. Mynhardt, A. Roux, and L. E. Teshima. Connected
$k$-dominating graphs. \emph{Discrete Math.}, 342:145--151, 2019.

\bibitem {26CM-MT}C. M. Mynhardt and L. E. Teshima. A note on some variations
of the $\gamma$-graph. \emph{J. Combin. Math. Combin. Comput.}, 104:217--230, 2018.

\bibitem {26CM-Nishimura}N. Nishimura. Introduction to reconfiguration.
\emph{Algorithms} (Basel, Paper No. 52), 11, 25 pp., 2018.

\bibitem {26CM-Sub}N. S. Sridharan and K. Subramanian. $\gamma$-Graph of a
graph. \emph{Bull. Kerala Math. Assoc}., 5:17--34, 2008.

\bibitem {26CM-Srid}N. Sridharan and K. Subramanian. Trees and unicyclic
graphs are $\gamma$-graphs. \emph{J. Combin. Math. Combin. Comput.},
69:231-236, 2009.

\bibitem {26CM-SMN}A. Suzuki, A. E. Mouawad, and N. Nishimura. Reconfiguration
of dominating sets. In \emph{COCOON} (Z. Cai, A. Zelikovsky, and A. Bourgeois,
Eds.), \emph{LNCS}, 8591:405--416, Springer, Heidelberg, 2014.

\bibitem {26CM-Tebbal}Y. Tebbal. \emph{On the Complexity of Reconfiguration of
Clique, Cluster Vertex Deletion, and Dominating Set}. Master's Thesis,
University of Waterloo, Waterloo, ON, Canada, 2015.

\bibitem {26CM-Vai}K. Vaidyanathan. \emph{Refiguring Graph Colorings}.
Master's Thesis, University of Waterloo, Waterloo, Canada, 2017.

\bibitem {26CM-West}D. West. \emph{Introduction to Graph Theory}, 2nd ed.
Prentice Hall, Upper Saddle River, NJ, 2001.

\bibitem {26CM-WrochnaDiss}M. Wrochna. \emph{The topology of solution spaces
of combinatorial problems}. Doctoral dissertation, University of Warsaw, 2018.
\end{thebibliography}
\end{document}